\documentclass[a4paper,dvipsnames,fleqn]{article}

\usepackage{amsmath}
\usepackage{amssymb}
\usepackage{amsthm}
\usepackage[utf8]{inputenc}
\allowdisplaybreaks
\usepackage{bm}
\usepackage{enumerate}
\usepackage[backend=biber,bibstyle=alphabetic,citestyle=alphabetic,sorting=nyt,maxnames=6,giveninits=false]{biblatex}
	\AtEveryBibitem{\clearfield{doi}}
	\AtEveryBibitem{\clearfield{isbn}}
	\AtEveryBibitem{\clearfield{issn}}
	\AtEveryBibitem{\clearfield{url}}
	\renewbibmacro{in:}{\ifentrytype{article}{}{\printtext{\bibstring{in}\intitlepunct}}}
	\AtEveryBibitem{\clearfield{number}}
	\AtEveryBibitem{\ifentrytype{book}{\clearfield{pages}}{}}
	\DeclareFieldFormat[article,incollection,inbook,unpublished]{title}{#1}
	\DeclareFieldFormat[article]{volume}{\textrm{#1}}
	\DeclareFieldFormat[article]{pages}{#1}
	\DeclareFieldFormat[book]{edition}{#1}
\usepackage[british]{babel}
\usepackage{csquotes}
\usepackage[UKenglish]{isodate}
\usepackage{dsfont}
\usepackage{enumitem}
\usepackage[OT2,OT1]{fontenc}
\usepackage[cal=cm,scr=boondoxo]{mathalfa}
\usepackage{pifont}
\usepackage{stmaryrd}
\usepackage{textcomp}
\usepackage{tikz}
	\usetikzlibrary{arrows}
	\usetikzlibrary{patterns}
\usepackage{tikz-cd}
\usepackage[textwidth=3cm,colorinlistoftodos]{todonotes}
\usepackage{xfrac}
\usepackage[colorlinks=true,allcolors=purple,pdfencoding=auto, psdextra]{hyperref}
\usepackage{cleveref}
	\crefformat{equation}{#2(#1)#3}
\numberwithin{equation}{section}			% Enumeration of equations
\swapnumbers						% Enumeration of theorem-like environments

\newcommand\cyr{%					% cyrillic font
\renewcommand\rmdefault{wncyr}%
\renewcommand\sfdefault{wncyss}%
\renewcommand\encodingdefault{OT2}%
\normalfont
\selectfont}
\DeclareTextFontCommand{\textcyr}{\cyr}

\newcounter{Enum}				% Enumerated list
\newenvironment{Enumerate}{\begin{enumerate}[label={\rm({\roman*})}]}{\end{enumerate}}

\newcommand{\descriptionlabelsave}{}		% Itemized list

\newcounter{StepsCount}				% Enumerated list with no indentation (e.g. steps in proof)

\newcounter{StepsRefCount}

\theoremstyle{plain}
	\newtheorem{lemma}{Lemma}[section]
	\newtheorem{proposition}[lemma]{Proposition}
	\newtheorem{theorem}[lemma]{Theorem}
	\newtheorem{corollary}[lemma]{Corollary}
	\newcommand{\GenericTheoremName}{}\newtheorem{generictheorem}[lemma]{\GenericTheoremName}
\theoremstyle{definition}
	\newtheorem{definition}[lemma]{Definition}
	\newcommand{\GenericDefinitionName}{}\newtheorem{genericdefinition}[lemma]{\GenericDefinitionName}
	\newtheorem{notation}[lemma]{Notation}
\theoremstyle{remark}
	\newtheorem{remark}[lemma]{Remark}
	\newtheorem{example}[lemma]{Example}
	\newcommand{\GenericRemarkName}{}\newtheorem{genericremark}[lemma]{\GenericRemarkName}
\newcommand{\mc}[1]{{\mathcal{#1}}}			% --- abbreviation ---
\newcommand{\ms}[1]{{\mathscr{#1}}}			% --- abbreviation ---
\newcommand{\bb}[1]{{\mathbb{#1}}}			% --- abbreviation ---
\DeclareMathOperator{\RE}{Re}				% real part

\DeclareMathOperator{\IM}{Im}				% imaginary part

\DeclareMathOperator{\supp}{supp}

\DeclareMathOperator{\tr}{tr}

\DeclareMathOperator{\sgn}{sgn}

\newcommand{\Side}[1]{\hfill{#1}\kern10pt}		% text put on the right side of line with offset

\newcommand{\smmatrix}[4]{\Bigl(			% small matrix for use in textline
\begin{smallmatrix}
\hspace*{-0.2ex} #1 \hspace*{0.2ex} & \hspace*{0.2ex} #2 \hspace*{-0.2ex}
\\[0.5ex]
\hspace*{-0.2ex} #3 \hspace*{0.2ex} & \hspace*{0.2ex} #4 \hspace*{-0.2ex}
\end{smallmatrix}
\Bigr)}
\newcommand{\Dummy}{\text{\textvisiblespace\kern1pt}}	% Platzhaltersymbol fuer Funktionsargumente
\newcommand{\Smallo}{{\rm o}}				% small o
\newcommand{\BigO}{{\rm O}}				% big o

\newcommand{\DS}{\mid\mkern3mu}				% delimiter for set definition
\newcommand{\DP}{{.\kern5pt}}				% delimiter for predicate formula
\newcommand{\DF}{\colon}				% delimiter for function domain/codomain
\newcommand{\DE}{\mathrel{\mathop:}=}			% defining equality
\newcommand{\ED}{=\mathrel{\mathop:}}			% defining equality
\newcommand{\DD}{\mkern4mu\mathrm{d}}			% distance and rm-d for integration differential

\begin{document}

\begin{flushleft}
	{\Large\bf Growth estimates for Nevanlinna matrices of \\[1ex]
	order larger than one half}
	\\[5mm]
	\textsc{
	Jakob Reiffenstein
	\hspace*{-14pt}
		\renewcommand{\thefootnote}{\fnsymbol{footnote}}
		\setcounter{footnote}{2}
			\footnote{
	The author was supported by the project I 4600 of the Austrian
Science Fund (FWF), and by the Sverker Lerheden foundation.}
\setcounter{footnote}{0}
	} \\[1ex]
	{\small
	\textbf{Abstract.}
Our objects of study are two-dimensional canonical systems that arise from indeterminate Hamburger moment problems and associated half-line Jacobi operators in limit circle case. The monodromy matrix of such a system coincides, up to a permutation of its entries, with the Nevanlinna matrix of the associated moment problem. Its growth relates to the density of eigenvalues of self-adjoint realisations of the system and the Jacobi operator, respectively. \\[1ex]
The order of the Nevanlinna matrix is known to be at most 1. In the case of ``large'' order, meaning order greater than one half, determining the growth of the monodromy matrix is known to be harder than for ``small'' order, i.e., order less than one half. As our main result, we establish an explicit new lower bound for data featuring a certain kind of monotonicity, which correctly describes the growth in the case of large order. 
%It shows that the same condition that is sufficient for determining the growth in the case of small order, falls short of doing so for large order. 
Moreover, we compute the order of the Nevanlinna matrix of a limit circle Jacobi matrix with two-term power asymptotics in a critical case.
	\\[3mm]
	}
\end{flushleft}

%%%%%%%%%%%%%%%%%%%%%%%%%%%%%%%%%%%%%%%%%%%%%%%%%%%%%%%%%%%%%%%%%%%%%%%%%%%%%%%%%%%%%%%%%
% disable in final version
%\pagenumbering{roman}
%\fbox{
%\parbox{100mm}{
%\hspace*{0pt}\\[1mm]
%\centerline{{\Large\ding{45}}\quad\,{\large\sc Draft}%\quad{\Large\ding{45}}}
%\hspace*{0pt}\\[-2mm]
%\textcircledP\ \ Preliminary version Fri 8 Jul 2023 16:46
%\\[2mm]
%\cleanlookdateon
%\hspace*{5mm} Compilation date: \today
%\\[2mm]
%\ding{233}\quad Use pdflatex/biber to compile
%\\[-1mm]
%}
%}
%\tableofcontents
%\listoftodos
%\newpage
%\pagenumbering{arabic}
%\setcounter{page}{1}
%
%%%%%%%%%%%%%%%%%%%%%%%%%%%%%%%%%%%%%%%%%%%%%%%%%%%%%%%%%%%%%%%%%%%%%%%%%%%%%%%%%%%%%%%%%

%---------
%   TEXTBODY
%---------

%\newpage

%
%
%
%\section[\textcolor{ForestGreen}{}]{}
%\section[\textcolor{Dandelion}{}]{}
%\section[\textcolor{BrickRed}{...}]{...}

\section[{Introduction}]{Introduction}

In the present paper, we consider two-dimensional \emph{canonical systems}
\begin{align}
\label{B36}
y'(t)=zJH(t)y(t), \qquad t \in (0,L) \text{ a.e.},
\end{align}
where $L \in (0,\infty )$, $J:=\smmatrix 0{-1}10$, $z \in \mathbb{C}$, and the \emph{Hamiltonian} $H$ is a $2 \times 2$ matrix-valued function on $(0,L)$ that satisfies 
\[
H(t) \geq 0, \quad \tr H(t)=1, \qquad t \in (0,L) \text{ a.e.,}
\]
and is of a particular discrete form (the precise definition is given further below). Following \cite{kac:2002} we call Hamiltonians of this kind \emph{Hamburger Hamiltonians}. This term is reminiscent of the fact that they form an equivalent model of Hamburger moment problems: the orthogonal polynomials $p_n$ of a moment sequence $(s_n)_{n=0}^\infty$ satisfy a three-term recurrence relation
\begin{align}
\label{K34}
b_n p_{n+1}(z) + a_n p_n(z)+b_{n-1}p_{n-1}(z) = zp_n (z)
\end{align}
with certain parameters $b_n>0$ and $a_n \in \mathbb{R}$. There exists a Hamburger Hamiltonian $H$ such that the minimal Jacobi operator with parameters $a_n,b_n$ is unitarily equivalent to the symmetric operator in $L^2(H)$ associated with \eqref{B36} and an initial condition at $0$\hspace{.9pt}\footnote{The initial condition can be chosen at will, but the Hamiltonian is different for different initial conditions.}. This process can be reversed, i.e., given a Hamburger Hamiltonian $H$ and an initial condition at $0$ there are unique Jacobi parameters $a_n,b_n$ such that the above mentioned operators are unitarily equivalent.

Properties of the moment sequence $(s_n)_{n=0}^\infty$ can be reformulated in terms of the corresponding Hamburger Hamiltonian $H$. Let us mention two prominent instances that are important for the present work:
\begin{itemize}
\item[$\rhd$] $(s_n)_{n=0}^\infty$ is indeterminate if and only if $L<\infty$;
\item[$\rhd$] If $(s_n)_{n=0}^\infty$ is indeterminate, its Nevanlinna matrix $W=\smmatrix ACBD$ and the monodromy matrix $W_H$ are related by $W_H=\smmatrix {C}{A}{-D}{-B}$.
\end{itemize}
We deal with the indeterminate case, and investigate specifically growth and zero distribution of the entries of the Nevanlinna matrix. Since these coincide, up to a possible minus sign, with the entries of the monodromy matrix of the system \eqref{B36}, their zero distribution has operator theoretic meaning: the set of zeroes of, say, $-B$ coincides with the spectrum of the self-adjoint operator arising from \eqref{B36} and the boundary conditions $(1,0)y(0)=(0,1)y(L)=0$. Hence, obtaining more detailed information about $A,B,C,D$ will yield knowledge about the spectra of self-adjoint extensions of the minimal Jacobi operator. As a simplest example: the (common) order of $A,B,C,D$ coincides with the (common) convergence exponent of the spectrum of all self-adjoint extensions.

A classical theorem of M. Riesz \cite{riesz:1923a}, or, alternatively, the Krein-de~Branges formula \cite{krein:1951,debranges:1961}, tells us that $W$ is of minimal exponential type. This implies that the order of $W$ is at most $1$, but its precise value is not given by these classical results. An abstract approach is given in \cite[Theorem~3.1]{berg.szwarc:2014}, equating the order of $W$ to the order of an entire function built from the coefficients of $p_n$, but this result is very hard (if not nearly impossible) to apply in practice.

There is a line of research, starting already in the early days of the study of moment problems, that deals with finding explicit upper and lower bounds for the order; a few examples are \cite{livshits:1939,berezanskii:1956,berg.szwarc:2014,romanov:2017,pruckner:blubb,pruckner.reiffenstein.woracek:sinqB-arXiv}. While a variety of upper bounds is available, there is essentially just one previously known lower bound, which stems from the comparison of the entries of $W$ to some entire functions derived from $(s_n)_{n=0}^\infty$ or coefficients of $p_n$ \cite{livshits:1939,berg.szwarc:2014}. In fact, in order to get an explicit estimate, all but the leading coefficients of $p_n$ need to be neglected, accounting for a very rough estimate. Nonetheless, under a regularity condition, this lower bound correctly describes the growth of $W$ if the order of $W$ is less than one half \cite[Theorem 5.3]{pruckner.reiffenstein.woracek:sinqB-arXiv}. For orders larger than one half this is, in general, not true anymore.

Our main aim in this paper is to derive and exploit a new lower bound for $\max_{|z|=r} \log  \|W(z)\|$. It is valid under a certain monotonicity assumption which, in particular, is fulfilled for Jacobi parameters having power asymptotics. If the order of $W$ is larger than one half, this lower bound coincides up to a multiplicative constant with a previously known upper bound. Hence we determine the exact growth of $\max_{|z|=r} \log  \|W(z)\|$ in this case, and consequently we can compute the order.
 We present the new lower bound and its consequences in the canonical systems setting as well as in the Jacobi matrix setting. 
 
 We point out that familiarity with canonical systems is not required to read the results in terms of Jacobi parameters, even though the proofs make use of recent work on canonical systems \cite{langer.reiffenstein.woracek:kacest-arXiv}.
\newline

We present our main results on canonical systems in Subsection \ref{SS11}. \Cref{K32} gives a lower bound for the monodromy matrix of any Hamburger Hamiltonian which ``rotates'' in a monotone way. \Cref{K29} combines this with previously known upper and lower bounds to give an overview of possible behaviours of $W$ when its order is larger than one half, emphasising the influence of the ``direction of rotation'' of $H$.
The proof of \Cref{K32} is given in \Cref{S3} along with \Cref{K61}, where an exceptional case is studied.

In Subsection \ref{SS12} we state our results on Jacobi matrices. We study parameters $a_n,b_n$ satisfying two-term power asymptotics in the critical case, i.e., $\lim_{n \to \infty} \frac{a_n}{b_n}=\pm 2$. For such parameters an explicit characterisation of the occurrence of limit circle or limit point case is available \cite{pruckner:blubb} (a third term in the expansion is needed in the doubly critical case). However, the precise order of the Nevanlinna matrix remained unknown for orders larger than one half. Using our results for canonical systems we determine the order of $W$ in this case, cf. \Cref{K10}. An overview of all cases is given in \Cref{K9}. In \Cref{S4} we consider more general Jacobi parameters: instead of power asymptotics we use the notion of $\gamma$-temperedness as in \cite{swiderski.trojan:2023}. We determine $\max_{|z|=r} \log  \|W(z)\|$ for such parameters, cf. \Cref{K78}, and deduce \Cref{K10} from this. \Cref{K78} is proved by passing over to the canonical systems setting, which is a rather delicate and technical task.

\subsection*{Hamburger Hamiltonians}
Let $(l_j)_{j=1}^\infty$ be a summable sequence of positive numbers and $(\phi_j)_{j=1}^\infty$ be a sequence of real
	numbers satisfying $\phi_{j+1}-\phi_j \not\equiv 0 \mod \pi$ for all $j \in \bb N$. Using the notation
	\begin{align*}
	x_0 &\DE 0, \qquad x_n \DE\sum_{j=1}^n l_j, \qquad n \in \mathbb{N} , \\  L &\DE\sum_{j=1}^\infty l_j <\infty
	\end{align*}
	and
\[
\xi_\phi \DE \binom{\cos \phi}{\sin \phi},
\]
	we define a Hamiltonian on $(0,L)$ by setting 
	\begin{align}
		H_{l,\phi}(t)\DE\xi_{\phi_j}\xi_{\phi_j}^\top\quad\text{for }j\in\bb N\text{ and }
		x_{j-1} \leq t< x_j
		.
	\end{align}
	This can be illustrated by
	\begin{center}
	\begin{tikzpicture}[x=1.2pt,y=1.2pt,scale=0.8,font=\fontsize{8}{8}]
		\draw[thick] (10,30)--(215,30);
		\draw[dotted, thick] (215,30)--(270,30);
		\draw[thick] (10,25)--(10,35);
		\draw[thick] (70,25)--(70,35);
		\draw[thick] (120,25)--(120,35);
		\draw[thick] (160,25)--(160,35);
		\draw[thick] (190,25)--(190,35);
		\draw[thick] (210,25)--(210,35);
		\draw[thick] (270,25)--(270,35);
		\draw (40,44) node {${\displaystyle \xi_{\phi_1}\xi_{\phi_1}^\top}$};
		\draw (95,44) node {${\displaystyle \xi_{\phi_2}\xi_{\phi_2}^\top}$};
		\draw (140,44) node {${\displaystyle \xi_{\phi_3}\xi_{\phi_3}^\top}$};
		\draw (177,43) node {${\cdots}$};
		\draw (-20,30) node {\large $H_{l,\phi}\!:$};
		\draw (10,18) node {${\displaystyle x_0}$};
		\draw[dashed,stealth-stealth] (11,26)--(69,26);
		\draw (40,21) node {${l_1}$};
		\draw (70,18) node {${x_1}$};
		\draw[dashed,stealth-stealth] (71,26)--(119,26);
		\draw (95,21) node {${l_2}$};
		\draw (120,18) node {${x_2}$};
		\draw[dashed,stealth-stealth] (121,26)--(159,26);
		\draw (140,21) node {${l_3}$};
		\draw (160,18) node {${x_3}$};
		\draw (195,18) node {${\cdots}$};
		\draw (270,18) node {${\displaystyle L}<\infty$};
	\end{tikzpicture}
	\end{center}
	
	The numbers $l_j$ and $\phi_j$ are referred to as \emph{lengths} and \emph{angles} of the Hamiltonian. Note that the lengths are unique for each Hamburger Hamiltonian, while the angles are not: adding an integer multiple of $\pi$ to any $\phi_j$ leaves the Hamiltonian invariant. 
This aspect becomes relevant already in \Cref{K32}, which features a condition that might be satisfied for some, but not all sequences $(\phi_j)_{j=1}^\infty$ representing the same Hamiltonian.

As was mentioned earlier, the condition $L<\infty$ means that $H_{l,\phi}$ is integrable on $(0,L)$, and thus corresponds to an indeterminate moment problem. In the terminology of canonical systems, this means that Weyl's limit circle case takes place at $L$.
\newline

%Hence, when we make an assumption on ``the sequence of angles'' of a Hamburger Hamiltonian $H$, what we mean is that \emph{there exists} $(\phi_j)_{j=1}^\infty$ such that $H=H_{l,\phi}$ and that satisfies the .

%Given a Hamburger Hamiltonian $H_{l,\phi}$, we will sometimes ask for $(\phi_j)_{j=1}^\infty$ to satisfy some condition $(\ast )$. We point out that even if $(\ast )$ does not hold for $(\phi_j)_{j=1}^\infty$, it might be satisfied by another sequence $(\tilde \phi_j)_{j=1}^\infty$ for which $H_{l,\phi}=H_{l,\tilde \phi}$.

%This means that 
%In order to avoid confusion, we point out that for a given Hamburger Hamiltonian $H_{l,\phi}$ and a given condition $(\ast )$ to be satisfied
	
\noindent Let $W_H(t;z)$ be the fundamental solution of (the transpose of) \eqref{B36}:
\begin{align*}
\begin{cases}
\frac{d}{dt} W_H(t;z)J=zW_H(t;z)H(t), &t \in [0,L], \\
W_H(0;z)=I.
\end{cases}
\end{align*}
Since limit circle case takes place at $L$, the solution $W_H$ exists up to $t=L$, and with some abuse of notation we set
\[
W_H(z)\DE W_H(L;z), \qquad \rho_H \DE \limsup_{r \to \infty} \frac{\log \log \max_{|z|=r} \|W_H(z)\|}{\log r}.
\]
The function $W_H=(w_{H;ij})_{i,j=1}^2$ is called the \emph{monodromy matrix}. It is an entire function of exponential type zero, and its \emph{order} $\rho_H$ thus satisfies $0\leq \rho_H \leq 1$. 

\subsection*{Auxiliary notions}
\begin{enumerate}
\item A function $\ms a\DF[1,\infty)\to(0,\infty)$ is called \emph{regularly varying} at $\infty$, if it is measurable and satisfies
	\[
		\forall \lambda\in(0,\infty)\DP \lim_{r\to\infty}\frac{\ms a(\lambda r)}{\ms a(r)}=\lambda^\alpha		
	\]
for some $\alpha\in\bb R$ called the \emph{index} of $\ms a$.
% We write $\ind \ms a \DE \alpha$. 
\item Assume that $\ms a$ is regularly varying with positive index. Then
\begin{align}
\label{K67}
\ms a^-(x) \DE \sup \big\{t\in[1,\infty)\DS \ms a(t)<x\big\}
\end{align}
is an asymptotic inverse of $\ms a$, i.e.,
\begin{equation}
	\label{B17}
		(\ms a\circ\ms a^-)(x)\sim(\ms a^-\circ\ms a)(x)\sim x
		.
	\end{equation}
	\item A sequence $(c_n)_{n=1}^\infty$ of positive numbers is called regularly varying if
	\[
		\forall \lambda\in(0,\infty)\DP \lim_{n\to\infty}\frac{c_{\lfloor \lambda n\rfloor}}{c_n}=\lambda^\alpha		
	\]
for some $\alpha\in\bb R$. We call $\alpha$ the \emph{index} of $(c_n)_{n=1}^\infty$.
\item Consider functions $f,g: X \to (0,\infty)$, where $X$ is any set.
We use the following notation to compare $f$ and $g$.
\begin{itemize}
\item[$\rhd$] We write $f \lesssim g$ (or $f(x) \lesssim g(x)$) if there exists $C>0$ such that $f(x) \leq Cg(x)$ for all $x \in X$. By $f \gtrsim g$ we mean $g \lesssim f$, and $f \asymp g$ is used for ``$f \lesssim g$ and $f \gtrsim g$''.
\item[$\rhd$] Let $Y \subseteq X$. We write ``$f \lesssim g$ on $Y$'' if $f|_Y \lesssim g|_Y$.
\item[$\rhd$] If $X=[1,\infty)$ or $X=\bb N$ (or any other directed set) we write ``$f(x) \lesssim g(x)$ for sufficiently large $x$'' to say that there exists $x_0 \in X$ such that $f \lesssim g$ 
on $\{x \in X \DF \, x \geq x_0 \}$. 
\item[$\rhd$] The notation $f \sim g$ stands for $\lim_{x \to \infty} \frac{f(x)}{g(x)}=1$. We write $f=\BigO(g)$ for $\limsup_{x \to \infty} \frac{f(x)}{g(x)}<\infty$, while $f=\Smallo (g)$ means $\lim_{x \to \infty} \frac{f(x)}{g(x)}=0$.
\end{itemize}

\end{enumerate}

\subsection[{Main results in terms of Hamiltonian parameters}]{Main results in terms of Hamiltonian parameters}
\label{SS11}
Our central assumption on a Hamburger Hamiltonian $H$ is that is lengths and the step sizes of its angles are comparable to regularly varying functions $\ms d_l, \ms d_\phi$ (in the sense that \eqref{K30} holds).
 It is known from \cite[Theorem 5.3] {pruckner.reiffenstein.woracek:sinqB-arXiv} that $\max_{|z|=r}\log\|W_H(z)\|$ is then determined up to $\asymp$ by the inverse function of $\frac{1}{\ms d_l \ms d_{\phi}}$, provided that $\ms d_l$ and $\ms d_\phi$ decay fast enough, i.e., the sum of their indices of regular variation is less than $-2$. In contrast to this, if the sum of these indices is larger than $-2$, the growth of $W_H$ depends also on the signs of the angle increments. Our main result, \Cref{K32}, gives a new lower bound in the case of constant sign of angle increments, i.e., monotone angles. 

\begin{theorem}
\label{K32}
Consider a Hamburger Hamiltonian $H$ in limit circle case. Let $\ms d_l,\ms d_\phi: [1,\infty) \to (0,\infty)$ be regularly varying with indices $-\delta_l$ and $-\delta_\phi$, respectively. Assume that $\delta_\phi \in (0,1)$, and suppose that angles $(\phi_j)_{j=1}^\infty$ of $H$ can be chosen such that
\begin{alignat}{2}
\label{K95}
&\rhd \quad && l_j \gtrsim \ms d_l(j), \quad |\phi_{j+1}-\phi_j| \asymp \ms d_\phi (j) \quad \text{for suff. large } j; \\[.5ex]
\nonumber
&\rhd \quad && \text{$(\phi_j)_{j=1}^\infty$ is eventually monotone.}
\end{alignat}
Then
\begin{align}
\label{K21}
\log |w_{H,22}(ir)| \gtrsim r\int_{\big[\frac{\ms d_\phi}{\ms d_l} \big]^-(r)}^\infty \ms d_l (x) \DD x.
\end{align}
for $r$ sufficiently large. In particular, $\rho_H \geq \frac{1-\delta_\phi}{\delta_l-\delta_\phi}$.
\end{theorem}

Note that $\delta_l$ is not less than $1$, due to $L<\infty$.  \newline

\Cref{K32} is most relevant when $1<\delta_l+\delta_\phi<2$ and $l_j \asymp \ms d_l(j)$, where also $\lesssim$ holds on the left hand side of \eqref{K95}. By \cite[Corollary 4.7]{pruckner.reiffenstein.woracek:sinqB-arXiv}, the lower bound from \eqref{K21} is then an upper bound, too. 

If $\delta_l+\delta_\phi>2$, the previously known lower bound, cf. \cite[Corollary 5.2]{pruckner.reiffenstein.woracek:sinqB-arXiv}, is better than the one from \Cref{K32}. We arrive at the following overview. 

\begin{theorem}
\label{K29}
Let $H=H_{l,\phi}$ be a Hamburger Hamiltonian in limit circle case, and suppose that 
\begin{align}
\label{K30}
l_j \asymp \ms d_l(j), \quad |\sin(\phi_{j+1}-\phi_j)| \asymp \ms d_\phi (j)
\end{align}
for sufficiently large $j$,
where $\ms d_l$ is regularly varying with index $-\delta_l$, and $\ms d_\phi$ is regularly varying with index $-\delta_\phi$. Assume that $\ms d_\phi(t)$ is $\sim$ to a nonincreasing function as $t \to \infty$. \\

Set 
\begin{align}
\label{K31}
\ms k(r) \DE \Big[\frac{1}{\ms d_l \ms d_\phi} \Big]^- (r), \qquad \ms m(r) \DE r\int_{\big[\frac{\ms d_\phi}{\ms d_l} \big]^-(r)}^\infty \ms d_l (x) \DD x.
\end{align}
Then, if
\[
1<\delta_l+\delta_\phi<2,
\]
the following statements hold.
	\begin{Enumerate}
	\item In general, for sufficiently large $r$ we have
	\[
	\ms k (r) \lesssim \max_{|z|=r}\log\|W_H(z)\| \lesssim \ms m(r).
	\]
	In particular, $\frac{1}{\delta_l+\delta_\phi} \leq \rho_H \leq \frac{1-\delta_\phi}{\delta_l-\delta_\phi}$.
	\item If there is $\psi \in \mathbb{R}$ such that $|\sin (\phi_j-\psi )| \lesssim |\sin (\phi_{j+1}-\phi_j )|$, and if $\delta_l>1$, then for sufficiently large $r$ we have
	\[
	\max_{|z|=r}\log\|W_H(z)\| \asymp \ms k (r).
	\]
	In particular, $\rho_H=\frac{1}{\delta_l+\delta_\phi}$.
	\item If $\phi_j$ is monotone with $|\phi_{j+1}-\phi_j| \asymp |\sin(\phi_{j+1}-\phi_j)| \asymp \ms d_\phi(j)$ for sufficiently large $j$, and if $\delta_\phi>0$, then for sufficiently large $r$ we have
	\[
	\max_{|z|=r}\log\|W_H(z)\| \asymp \ms m(r)
	\]
	In particular, $\rho_H=\frac{1-\delta_\phi}{\delta_l-\delta_\phi}$.
	\end{Enumerate}
\end{theorem}

\begin{remark}
\label{K37}
We expect that, in the situation of \Cref{K29}, any growth between $\ms k(r)$ and $\ms m(r)$ is possible, meaning that for any regularly varying function $\ms f$ with $\ms k \lesssim \ms f \lesssim \ms m$ there is a Hamburger Hamiltonian satisfying \eqref{K30}, for which
\[
	\max_{|z|=r}\log\|W_H(z)\| \asymp \ms f(r).
\]
Our idea how this may be achieved is to fix 
\[
l_j \DE \ms d_l(j), \quad |\phi_{j+1}-\phi_j| \DE \ms d_{\phi}(j)
\]
and assign signs to $\phi_{j+1}-\phi_j$ such that $\phi_j$ approximates $\ms u(j)$ for some smooth target function $\ms u$:
\[
\phi_1 \DE \frac{\pi}{2}, \qquad \phi_{j+1} \DE 
\begin{cases}
\phi_j+\ms d_{\phi}(j), & \text{if } \phi_j+\ms d_{\phi}(j)<\ms u(j), \\
\phi_j-\ms d_{\phi}(j), & \text{if } \phi_j+\ms d_{\phi}(j) \geq \ms u(j).
\end{cases}
\]
We predict that, if $\ms u$ and $\ms u'$ are regularly varying and the index of $\ms u$ is between $0$ and $1-\delta_{\phi}$, the estimate
\begin{align}
\label{K35}
\max_{|z|=r}\log\|W_H(z)\| \asymp \max \bigg\{\ms k(r),\, r\int_{\ms v^-(r)}^\infty \ms d_l (x) \DD x \bigg\}
\end{align}
will hold, where $\ms v(x) \DE \frac{u'(x)}{\ms d_l(x)}$.

The major challenge in proving this is to establish the upper estimate, whereas the lower estimate seems to be in reach by generalising the proof of \Cref{K32}.
\end{remark}

\subsection[{Main results in terms of Jacobi parameters}]{Main results in terms of Jacobi parameters}
\label{SS12}

Consider a Jacobi matrix $\ms J$ with parameters satisfying
\begin{align}
\label{K53}
\begin{split}
b_n &=n^{\beta_1} \Big(x_0+\frac{x_1}{n}+\chi_n \Big), \\
a_n &=n^{\beta_2} \Big(y_0+\frac{y_1}{n}+\mu_n \Big),
\end{split}
\end{align}
with $x_0>0$, $y_0 \neq 0$ and $\chi_n,\mu_n$ are ``small''. For suitable notions of ``small'', the occurrence of limit point or limit circle case is characterised in \cite{pruckner:blubb} in terms of the parameters $\beta_1,\beta_2$ and $x_0,x_1,y_0,y_1$. How small $\chi_n,\mu_n$ need to be depends on which of the following cases takes place:
\begin{itemize}
\item[$\rhd$] Large diagonal: 
	\[
		\beta_2>\beta_1\quad\text{or}\quad\Big(\beta_1=\beta_2\wedge|y_0|>2x_0\Big)
		;
	\]
\item[$\rhd$] Small diagonal:
	\[
		\beta_2<\beta_1\quad\text{or}\quad\Big(\beta_1=\beta_2\wedge|y_0|<2x_0\Big)
		;
	\]
\item[$\rhd$] Simply critical case: 
\[
\beta_2=\beta_1 \ED \beta \quad \wedge \quad |y_0|=2x_0 \quad  \wedge \quad \beta \neq \sigma \DE \frac{2x_1}{x_0}-\frac{2y_1}{y_0};
\]
\item[$\rhd$] Doubly critical case: 
\[
\beta_2=\beta_1 \, (=\beta) \quad \wedge \quad |y_0|=2x_0 \quad \wedge \quad \beta = \sigma.
\]
\end{itemize}

In the large and small diagonal cases, it suffices to take $\chi_n,\mu_n=\BigO(n^{-(1+\epsilon)})$ (for fixed $\epsilon >0$). Then $\ms J$ is in limit point case if its diagonal is large. If its diagonal is small, it is in limit circle case if and only if $\beta_1>1$, and if so, its Nevanlinna matrix is of order $\frac{1}{\beta_1}$, cf. \cite[Theorem~1]{pruckner:blubb}. \\
The case we are interested in is the simply critical case; we use the assumption $\sum_{n=1}^\infty \sqrt{n} \big(|\chi_n|+|\mu_n| \big)<\infty$. It follows from \cite[Theorems~9.1 and 9.2]{swiderski.trojan:2023} that the limit circle case occurs if and only if $\frac 32<\beta<\sigma$. About the order $\rho$ of the Nevanlinna matrix it is known that $\rho=\frac{1}{\beta}$ when $2 \leq \beta<\sigma$ and that $\frac{1}{\beta} \leq \rho \leq \frac{1}{2(\beta-1)}$ when $\frac 32 <\beta<\min \{2,\sigma \}$, cf. \cite[Theorem~2]{pruckner:blubb}. By showing that the Hamburger Hamiltonian associated with $\ms J$ has monotone angles, \Cref{K29} lets us determine the order in the missing case.

\begin{theorem}
\label{K10}
Let $\ms J$ be a Jacobi matrix whose parameters satisfy \eqref{K53} with $\sum_{n=1}^\infty \sqrt{n} \big(|\chi_n|+|\mu_n| \big)<\infty$. Assume that the simply critical case takes place and that 
\begin{align*}
\frac 32 < \beta < \min\{2,\sigma \}.
\end{align*}
Then $\ms J$ is in limit circle case, and its Nevanlinna matrix $W$ satisfies
\begin{align}
\label{K65}
 \max_{|z|=r} \log \|W(z)\| \asymp r^{\frac{1}{2 (\beta -1)}}, \qquad r \text{ sufficiently large}
\end{align}
and its order is $\frac{1}{2(\beta-1)}$.
\end{theorem}

For the proof of this result we refer to \Cref{K90}.

%Note that $\ms J$ is still in limit circle case when $2 \leq \beta <\sigma$; in that case it is known (\cite[Theorem~2]{pruckner:blubb}, using slightly stronger assumptions) that the order of $W$ is equal to $\frac{1}{\beta}$. \\
In the doubly critical case the assumptions on $\chi_n,\mu_n$ we used above are not enough to distinguish limit circle and limit point case. However, after refining the assumptions we can paint the following complete picture, which combines \cite[Theorems~1,2]{pruckner:blubb} with \Cref{K10}.

\begin{theorem}
\label{K9}
Let $\ms J$ be a Jacobi matrix with parameters satisfying \eqref{K53} with $\chi_n=\frac{x_2}{n^2}+\BigO(n^{-(2+\epsilon)})$ and $\mu_n=\frac{y_2}{n^2}+\BigO(n^{-(2+\epsilon)})$. Then $\ms J$ is in limit circle case if and only if one of the following cases takes place.
\begin{Enumerate}
\item Small diagonal case with $\beta_1>1$;
\item Simply critical case with $\frac 32<\beta_1<\sigma$;
\item Doubly critical case with $2<\beta_1=\sigma<\frac 32+\eta$, where
\begin{align*}
\eta \DE 2 \bigg( \frac{2x_2}{x_0}-\frac{2y_2}{y_0}+\frac{y_1}{y_0}-\frac{2x_1y_1}{x_0y_0}+\frac{2y_1^2}{y_0^2} \bigg).
\end{align*}
\end{Enumerate}
If $\ms J$ is in limit circle case, the order $\rho$ of its Nevanlinna matrix is equal to $\frac{1}{\beta_1}$, except in the simply critical case when $\frac 32<\beta_1<2$, where $\rho$ is equal to $\frac{1}{2(\beta_1-1)}$.
\end{theorem}

\begin{remark}
The number $\frac{1}{\beta_1}$, which equals the order in most of these cases, is also equal to the \emph{convergence exponent} of $(b_n)_{n=0}^\infty$, i.e.,
\[
\frac{1}{\beta_1}=\inf \Big\{\alpha>0 \DS \sum_{n=0}^\infty\frac 1{b_n^\alpha}<\infty\Big\}.
\]
To the best of our knowledge, the Jacobi matrices considered in \Cref{K10} are the first known examples where limit circle case holds and the order of the Nevanlinna matrix can be computed, but is different from the convergence exponent of $(b_n)_{n=0}^\infty$.
\end{remark}

\section[{Regularly varying Hamiltonian parameters and growth of \boldmath{$W_H$}}]{Regularly varying Hamiltonian parameters and growth of \boldmath{$W_H$}}
\label{S3}

After some preparation we are going to prove \Cref{K32} in Subsection \ref{SS32}. In \Cref{K61} we give an interesting application of the proof method in the exceptional case where $\delta_l+\delta_\phi=2$.

\subsection[{A method of estimating \boldmath{$W_H$} from below}]{A method of estimating \boldmath{$W_H$} from below}
Let us recall some notions from \cite{langer.reiffenstein.woracek:kacest-arXiv}. First, we introduce the notation
\begin{align}
\Omega_H (s,t)=\int_s^t H(x) \, dx.
\end{align}
If the reference to $H$ is unambiguous, we write $\Omega$ instead of $\Omega_H$. \\[1ex]

The following lemma, which combines Proposition 5.9, Theorem 5.3, and Theorem 3.4 from \cite{langer.reiffenstein.woracek:kacest-arXiv}, is our central tool for bounding the growth of $W_H$ from below.

\begin{lemma}
\label{X67}
	Assume we have points $0 \leq s_0<s_1<\ldots <s_k \leq L$ and let $c,r>0$ such that
	\begin{equation}\label{X61}
		\det\Omega(s_{j-1},s_j) \ge \frac{c}{r^2}, \qquad j\in\{1,\ldots,k\}.
	\end{equation}
	Then we have
	\[
		\log |w_{H,22}(ir)| \geq c_0 \cdot \Big(
		k\log 2-\Big[\log r+\log\frac{2L}{\sqrt c}\Big]\Big)
	\]
where $c_0>0$ is a universal constant (i.e., independent of $H$).
\end{lemma}

\Cref{X67} is an important ingredient in the proof of \Cref{K32}. We will apply it by constructing, for each sufficiently large $r>0$, points $s_0^{(r)},s_1^{(r)},\ldots,s_{k(r)}^{(r)}$ such that \eqref{X61} is satisfied ($c$ is chosen the same for every $r$). Assuming that $k(r)$, viewed as a function of $r$, grows faster than $\log r$, this strategy gives us the lower bound
\[
		\log |w_{H,22}(ir)| \gtrsim k(r), \qquad r \text{ sufficiently large}.
\]

We will frequently use the fact that $\det \Omega_H$ can be rewritten in the following way:
\begin{align}
\label{KX54}
\det\Omega_H(x_m,x_n) =\frac 12 \sum_{j,k=m+1}^n l_jl_k \sin^2(\phi_j-\phi_k).
\end{align}
This formula is a special case of \cite[Lemma 6.3]{langer.reiffenstein.woracek:kacest-arXiv}.

\subsection[{Proof of Theorem 1.1}]{Proof of \protect{\Cref{K32}}}
\label{SS32}

We start with two auxiliary lemmata. The first one is folklore, but due to lack of reference we provide a proof.

\begin{lemma}
\label{K12}
Let $\ms a : [1,\infty) \to (0,\infty)$ be regularly varying and let $k>0$. Then there exists $M>0$ such that
\begin{align}
\label{K7}
\ms a(x+h) \asymp \ms a (x), \qquad x \in [M,\infty), \quad h \in [-k,k].
\end{align}
In particular,
\begin{align}
\label{K8}
\int_x^y \ms a(s) \DD s \asymp (y-x)\ms a(x) \asymp \frac{y-x}{y'-x'} \int_{x'}^{y'} \ms a(s) \DD s.
\end{align}
on $\Big\{(x,y,x',y') \DF M\leq x<y<x+k \, \wedge \, M \leq x'<y' \, \wedge \, [x',y'] \subseteq [x-k,x+k] \Big\}$. 
\end{lemma}
\begin{proof}
Let $K \DE \big[ \frac{1}{k+1},k+1 \big]$. Then for $x \geq k+1$ and $h \in [-k,k]$ we have $\frac{x+h}{x} \in K$. By \cite[Theorem 1.2.1]{bingham.goldie.teugels:1989}, we have $\lim_{x \to \infty} \frac{\ms a(\lambda x)}{\ms a(x)}=\lambda^\alpha$ uniformly for $\lambda \in K$. Consequently, there exists $M \geq k+1$ such that
\[
\bigg|\frac{\ms a(x+h)}{\ms a(x)} - \bigg(\frac{x+h}{x}\bigg)^\alpha \bigg| < \frac 1{2(k+1)^\alpha}, \quad x \geq M.
\]
Since $\big(\frac{x+h}{x}\big)^\alpha \geq \frac 1{(k+1)^{\alpha }}$, \eqref{K7} follows. Now \eqref{K8} holds because of
\[
\int_x^y \ms a(s) \DD s \asymp (y-x)\ms a(x) = \frac{y-x}{y'-x'} \cdot (y'-x') \ms a(x) \asymp \frac{y-x}{y'-x'}\int_{x'}^{y'} \ms a(s) \DD s.
\]
\end{proof}

The main intuition in the proof of \Cref{K32} is the idea that, for monotone angles, $|\sin (\phi_j-\phi_k)|$ should stay away from zero for many pairs of indices $j,k$. This will make $\det \Omega (x_m,x_n)$ large enough to choose sufficiently many points for \Cref{X67} to provide the desired lower bound. In the following lemma, we develop the necessary tool for quantifying the intuition about $|\sin (\phi_j-\phi_k)|$.

\begin{lemma}
\label{K15}
Let $H$ be a Hamburger Hamiltonian. Then, for $m,n \in \mathbb{N}$ with $m<n$ and a permutation $\sigma$ of $\{m+1,\ldots,n\}$,
\begin{align*}
\det \Omega(x_m,x_n) \geq \frac 18 \bigg(\sum_{j=m+1}^n \min \{l_j,l_{\sigma (j)}\} |\sin (\phi_j-\phi_{\sigma (j)})| \bigg)^2.
\end{align*}
\end{lemma}
\begin{proof}
We use the elementary fact $\sin^2(a+b) \leq 2 [\sin^2(a)+\sin^2(b)]$. By \eqref{KX54},
\begin{align*}
\det &\Omega (x_m,x_n) =\frac 12 \sum_{j,k=m+1}^n l_jl_k \sin^2(\phi_j-\phi_k) \\
&=\frac 14 \sum_{j,k=m+1}^n \Big[l_jl_k \sin^2(\phi_j-\phi_k)+l_kl_{\sigma (j)} \sin^2(\phi_k-\phi_{\sigma (j)}) \Big] \\
&\geq \frac 14 \sum_{j,k=m+1}^n l_k \min \{l_j,l_{\sigma (j)}\}\big[ \sin^2(\phi_j-\phi_k)+\sin^2(\phi_k-\phi_{\sigma (j)}) \big] \\
&\geq \frac 18 \bigg(\sum_{k=m+1}^n l_k \bigg) \bigg( \sum_{j=m+1}^n \min \{l_j,l_{\sigma (j)}\}\sin^2(\phi_j-\phi_{\sigma (j)}) \bigg) \\
&\geq \frac 18 \bigg( \sum_{j=m+1}^n \sqrt{l_j \min\{l_j,l_{\sigma (j)}\}} |\sin(\phi_j-\phi_{\sigma (j)})| \bigg)^2 \\
&\geq \frac 18 \bigg( \sum_{j=m+1}^n \min\{l_j,l_{\sigma (j)}\} |\sin(\phi_j-\phi_{\sigma (j)})| \bigg)^2.
\end{align*}
\end{proof}

Before we come to the proof of \Cref{K32}, we restate its formulation for the convenience of the reader. 

\begin{theorem}
Consider a Hamburger Hamiltonian $H$ in limit circle case. Let $\ms d_l,\ms d_\phi: [1,\infty) \to (0,\infty)$ be regularly varying with indices $-\delta_l$ and $-\delta_\phi$, respectively. Assume that $\delta_\phi \in (0,1)$, and suppose that angles $(\phi_j)_{j=1}^\infty$ of $H$ can be chosen such that
\begin{alignat*}{2}
&\rhd \quad && l_j \gtrsim \ms d_l(j), \quad |\phi_{j+1}-\phi_j| \asymp \ms d_\phi (j) \quad \text{for suff. large } j; \\[.5ex]
&\rhd \quad && \text{$(\phi_j)_{j=1}^\infty$ is eventually monotone.}
\end{alignat*}
Then
\begin{align*}
\log |w_{H,22}(ir)| \gtrsim r\int_{\big[\frac{\ms d_\phi}{\ms d_l} \big]^-(r)}^\infty \ms d_l (x) \DD x.
\end{align*}
for $r$ sufficiently large. In particular, $\rho_H \geq \frac{1-\delta_\phi}{\delta_l-\delta_\phi}$.
\end{theorem}

\begin{proof}
With \Cref{X67} in mind, for given $r>0$ we need to construct sufficiently many points $s_0^{(r)},s_1^{(r)},\ldots,s_{k(r)}^{(r)}$ such that $\det \Omega (s_{j-1}^{(r)},s_j^{(r)}) \geq \frac{c}{r^2}$ for all $j \in \{1,\ldots,k(r)\}$ and some $c$ independent of $r$. We will use \Cref{K15} in order to obtain a good lower estimate for $\det \Omega (x_m,x_n)$ whenever $m,n$ are large enough and not too close to each other. The main difficulty will be the choice of $\sigma$ in every application of \Cref{K15}, and this takes up a significant part of the proof.

We may assume that $\ms d_\phi$ is continuous: If not, replace it by the function obtained from linearly interpolating the integer values of $\ms d_\phi$. Let 
\[
\Phi (t) \DE \int_1^{t} \ms d_\phi (s) \DD s, \qquad t \geq 1.
\]
The function $\Phi : [1,\infty) \to [0,\infty)$ is increasing and regularly varying with index $1-\delta_\phi>0$. In particular, it is bijective. 
\vspace{8pt}
	
	{\it Step 1.} For fixed $h>0$, the mean value theorem yields
\begin{align*}
&\Phi^{-1}(x+h)-\Phi^{-1}(x)=h \big(\Phi^{-1}\big)'(\xi_x) = \frac{h}{\Phi' \big(\Phi^{-1}(\xi_x)\big)}
\end{align*}
for some $\xi_x \in [x,x+h]$. Since $\lim_{x \to \infty} \frac{\xi_x}{x}=1$ and because $\Phi^{-1}$ and $\Phi'$ are regularly varying, we have
\begin{align*}
&\Phi^{-1}(x+h)-\Phi^{-1}(x) \sim 
\frac{h}{\Phi' \big(\Phi^{-1}(x)\big)} = \frac{h}{\ms d_\phi \big(\Phi^{-1}(x)\big)}.
\end{align*}
By Karamata's Theorem \cite[Theorem 1.5.11]{bingham.goldie.teugels:1989}, $\Phi (t) \sim \frac{t \ms d_\phi(t)}{1-\delta_\phi}$, and hence
\begin{align}
\label{K6}
&\Phi^{-1}(x+h)-\Phi^{-1}(x) \sim 
\frac{h}{1-\delta_\phi} \cdot \frac{\Phi^{-1}(x)}{x}.
\end{align}
\vspace{8pt}
	
	{\it Step 2.} The right hand side of \eqref{K6} is regularly varying with index $\frac{\delta_\phi}{1-\delta_\phi}>0$, and thus tends to $\infty$. Fix numbers $u,w,u',w' \in (0,1)$ with $u<w$, $u'<w'$, and $w-u<w'-u'$. Consider, for $j \in \mathbb{N}$, the discrete intervals
\[
I_j^+\DE \Big\{k \in \mathbb{N} \,:\, \Big\lceil \Phi^{-1}\big(j+u \big) \Big\rceil \leq k \leq \Big\lfloor \Phi^{-1}\big(j+w \big)\Big) \Big\rfloor \Big\}
\]
and
\[
I_j^- \DE \Big\{l \in \mathbb{N} \,:\, \Big\lceil \Phi^{-1}\big(j-w' \big)\Big) \Big\rceil \leq l \leq \Big\lfloor \Phi^{-1}\big(j-u' \big)\Big) \Big\rfloor \Big\}.
\]
By \eqref{K6}, the cardinalities of $I_j^+$ and $I_j^-$ satisfy
\begin{align*}
\lim_{j \to \infty} \frac{|I_j^+|}{|I_j^-|} = \frac{w-u}{w'-u'} < 1.
\end{align*}
In particular, $|I_j^+|<|I_j^-|$ for all sufficiently large $j$.
\vspace{8pt}
	
	{\it Step 3.}
For the following estimates, we need a (finite) number of properties that only hold for sufficiently large indices. We assume implicitly that each variable that occurs is large enough, and we understand $f \asymp g$ as ``$f(x) \asymp g(x)$ for sufficiently large $x$''. 

Given sufficiently large $n,m \in \mathbb{N}$ with $n \geq m+2$, there exists a permutation $\sigma$ of $\big\{\big\lceil \Phi^{-1}\big( m\big) \big\rceil, \ldots, \big\lceil \Phi^{-1}\big(n\big) \big\rceil \big\}$ with the property $\sigma (I_j^+) \subseteq I_j^-$ for all $m+1 \leq j \leq n$. Our goal is to choose $u,w,u',w'$ such that
\begin{align}
\label{K14}
|\sin (\phi_k-\phi_{\sigma (k)})| \asymp 1, \quad m,n \text{ suff. large}, \quad m+1 \leq j \leq n, \quad k \in I_j^+.
\end{align} 
Let $k \in I_j^+$. Then
\begin{align}
\nonumber
u+&u' \leq \Phi \Big(\Big\lceil \Phi^{-1}\big(j+u \big)\Big) \Big\rceil \Big) - \Phi \Big(\Big\lfloor \Phi^{-1}\big(j-u' \big)\Big) \Big\rfloor \Big) \leq \Phi(k)-\Phi(\sigma (k)) \\[.5ex]
\label{K5}
&\leq \Phi \Big(\Big\lfloor \Phi^{-1}\big(j+w \big)\Big) \Big\rfloor \Big) - \Phi \Big(\Big\lceil \Phi^{-1}\big(j-w' \big)\Big) \Big\rceil \Big) \leq w+w'.
\end{align}
Our next goal is to obtain similar bounds for $|\phi_k-\phi_{\sigma (k)}|$.  By \Cref{K12}, 
\begin{align*}
\Phi (k)-\Phi (\sigma (k)) &=\int_{\sigma (k)}^{k} \ms d_\phi(s) \DD s = \sum_{i=\sigma (k)}^{k-1} \int_i^{i+1} \ms d_\phi (x) \DD x \asymp \sum_{i=\sigma (k)}^{k-1} \ms d_\phi (i) \\
&\asymp \sum_{i=\sigma (k)}^{k-1} \big|\phi_{i+1}-\phi_i \big| = \big|\phi_k-\phi_{\sigma (k)}\big|
\end{align*}
and hence there exist constants $c,c'>0$ such that
\begin{align*}
c \big( \Phi (k)-\Phi (\sigma (k))\big) \leq |\phi_k-\phi_{\sigma (k)}| \leq c' \big( \Phi (k)-\Phi (\sigma (k))\big)
\end{align*}
for $m,n$ large, $m+1 \leq j \leq n$ and $k \in I_j^+$. 
Remembering \eqref{K5}, we see that 
\begin{align}
\label{K13}
c (u+u') \leq \phi_k-\phi_{\sigma (k)} \leq c' (w+w')
\end{align}
for $m,n$ large, $m+1 \leq j \leq n$ and $k \in I_j^+$. 

Let $w\DE w' \DE \frac{\pi}{4} \min \{1, \frac{1}{c'} \}$ and $u=\frac w2$, $u'=\frac w4$. Due to \eqref{K13}, there is $\epsilon>0$ such that $\phi_k-\phi_{\sigma (k)} \in [\epsilon, \frac{\pi}{2}]$. Hence \eqref{K14} holds for this choice of $u,w,u',w'$.
\vspace{8pt}
	
	{\it Step 4.}
\noindent Recalling \eqref{K14},  \Cref{K15} leads to
\begin{align}
\label{K16}
\begin{split}
&\sqrt{\det \Omega \big(x_{\lceil \Phi^{-1}( m)\rceil},x_{\lceil \Phi^{-1}( n) \rceil} \big)} \\
&\gtrsim \sum_{k=\lceil \Phi^{-1}( m)\rceil+1}^{\lceil \Phi^{-1}( n)\rceil} \min \{l_k,l_{\sigma (k)} \} |\sin (\phi_k-\phi_{\sigma (k)})| \\
&\gtrsim \sum_{j=m+1}^{n} \sum_{k \in I_j^+} \min \{l_k,l_{\sigma (k)} \} \gtrsim \sum_{j=m+1}^{n} \sum_{k \in I_j^+} \min \{\ms d_l(k),\ms d_l(\sigma (k)) \}.
\end{split}
\end{align}
The function $\ms d_l$ is regularly varying with negative index. By \cite[Theorem 1.5.6 (iii)]{bingham.goldie.teugels:1989} and due to $\sigma(k)<k$ we have $\ms d_l (\sigma(k) \lesssim \ms d_l(k)$ for sufficiently large $k$. Using also \Cref{K12} we have 
\begin{align*}
&\sum_{j=m+1}^{n} \sum_{k \in I_j^+} \min \{\ms d_l(k),\ms d_l(\sigma (k)) \} \asymp \sum_{j=m+1}^{n} \sum_{k \in I_j^+} \ms d_l(k) \\
&\asymp \sum_{j=m+1}^{n} \sum_{k \in I_j^+} \int_k^{k+1} \ms d_l(x) \DD x = \sum_{j=m+1}^{n} \int_{\lceil \Phi^{-1}(j+u) \rceil}^{1+\lfloor \Phi^{-1}(j+w) \rfloor} \ms d_l(x) \DD x \\
&\asymp \sum_{j=m+1}^{n} \int_{\Phi^{-1}(j+u)}^{\Phi^{-1}(j+w)} \ms d_l(x) \DD x
\end{align*}
for $j$ sufficiently large. Combining this with \eqref{K16} yields
\begin{align}
\label{K17}
&\sqrt{\det \Omega \big(x_{\lceil \Phi^{-1}( m)\rceil},x_{\lceil \Phi^{-1}( n) \rceil} \big)} \gtrsim  \sum_{j=m+1}^{n} \int_{\Phi^{-1}(j+u)}^{\Phi^{-1}(j+w)} \ms d_l(x) \DD x
\end{align}
\vspace{8pt}
	
	{\it Step 5.} We further refine \eqref{K17}.
For all $\xi < \zeta$, substituting $x=\Phi^{-1}(t)$ yields
\begin{align}
\label{K19}
\int_{\Phi^{-1}(\xi)}^{\Phi^{-1}(\zeta)} \ms d_l(x) \DD x = \int_\xi^\zeta \frac{\ms d_l \big(\Phi^{-1}(t) \big)}{\Phi' \big(\Phi^{-1}(t) \big)} \DD t = \int_\xi^\zeta \frac{\ms d_l \big(\Phi^{-1}(t) \big)}{\ms d_\phi \big(\Phi^{-1}(t) \big)} \DD t.
\end{align}
Using \Cref{K12} again, this leads to
\begin{align}
\label{K18}
\begin{split}
&\sqrt{\det \Omega \big(x_{\lceil \Phi^{-1}( m)\rceil},x_{\lceil \Phi^{-1}( n) \rceil} \big)} \gtrsim \sum_{j=m+1}^{n} \int_{j+u}^{j+w} \frac{\ms d_l \big(\Phi^{-1}(t) \big)}{\ms d_\phi \big(\Phi^{-1}(t) \big)} \DD t \\
&\asymp \sum_{j=m+1}^{n} \frac{(j+w)-(j+u)}{(j+1)-j}\int_{j}^{j+1} \frac{\ms d_l \big(\Phi^{-1}(t) \big)}{\ms d_\phi \big(\Phi^{-1}(t) \big)} \DD t \\
&\asymp \int_{m+1}^{n+1} \frac{\ms d_l \big(\Phi^{-1}(t) \big)}{\ms d_\phi \big(\Phi^{-1}(t) \big)} \DD t \asymp \int_m^n \frac{\ms d_l \big(\Phi^{-1}(t) \big)}{\ms d_\phi \big(\Phi^{-1}(t) \big)} \DD t.
\end{split}
\end{align}

\vspace{8pt}
	
	{\it Step 6.}
For $r>0$ set
\[
m_0(r) \DE 1+\max \bigg\{m \in \mathbb{N} \DF \int_m^{m+1} \frac{\ms d_l \big(\Phi^{-1}(t) \big)}{\ms d_\phi \big(\Phi^{-1}(t) \big)} \DD t \geq \frac 1r \bigg\}
\]
and define, by induction,
\[
m_n(r) \DE \min \bigg\{m \in \mathbb{N} \DF m> m_{n-1}(r) \wedge \int_{m_{n-1}(r)}^m \frac{\ms d_l \big(\Phi^{-1}(t) \big)}{\ms d_\phi \big(\Phi^{-1}(t) \big)} \DD t \geq \frac 1r \bigg\}
\]
where $\min \emptyset \DE \infty$. For $n \geq 1$ with $m_n(r)<\infty$ we have, by \eqref{K18},
\begin{align}
\label{K22}
&\det \Omega \big(x_{\lceil \Phi^{-1}( m_{n-1}(r))\rceil},x_{\lceil \Phi^{-1}( m_n(r)) \rceil} \big) \gtrsim \frac 1{r^2}.
\end{align}
Let $k (r) \DE \min \{n \in \mathbb{N} \DF m_n(r)=\infty \}$. For large enough $r>0$, the ($k(r)$ many) points
\[
0, \, x_{\lceil \Phi^{-1}( m_1(r))\rceil},\, x_{\lceil \Phi^{-1}( m_2(r))\rceil}, \ldots , x_{\lceil \Phi^{-1}( m_{k(r)-2}(r))\rceil},\, L
\]
satisfy the condition of \Cref{X67}, and this yields that 
\begin{align}
\label{K68}
\log |w_{H,22}(ir)| \gtrsim  k (r)+ \BigO(\log r).
\end{align}
\vspace{8pt}
	
	{\it Step 7.}
We shall estimate $k(r)$ from below. By definition of $m_0(r)$, for $1 \leq n < k(r)$ it holds that
\begin{align}
\label{K54}
\frac 1r \leq \int_{m_{n-1}(r)}^{m_n(r)} \frac{\ms d_l \big(\Phi^{-1}(t) \big)}{\ms d_\phi \big(\Phi^{-1}(t) \big)} \DD t \leq \frac 1r + \int_{m_n(r)-1}^{m_n(r)} \frac{\ms d_l \big(\Phi^{-1}(t) \big)}{\ms d_\phi \big(\Phi^{-1}(t) \big)} \DD t \leq \frac 2r.
\end{align}
In addition we have
\[
\int_{m_{k(r)-1}(r)}^{\infty} \frac{\ms d_l \big(\Phi^{-1}(t) \big)}{\ms d_\phi \big(\Phi^{-1}(t) \big)} \DD t \leq \frac 1r.
\]
This leads to
\begin{align}
\label{K20}
\begin{split}
\frac{k (r)}{r} &\gtrsim \sum_{n=1}^{k (r)} \int_{m_{n-1}(r)}^{m_n(r)} \frac{\ms d_l \big(\Phi^{-1}(t) \big)}{\ms d_\phi \big(\Phi^{-1}(t) \big)} \DD t = \int_{m_0(r)}^\infty \frac{\ms d_l \big(\Phi^{-1}(t) \big)}{\ms d_\phi \big(\Phi^{-1}(t) \big)} \DD t.
\end{split}
\end{align}
By \Cref{K12},
\begin{align*}
\frac{\ms d_l \big(\Phi^{-1}(m_0 (r)) \big)}{\ms d_\phi \big(\Phi^{-1}(m_0(r)) \big)} \asymp \int_{m_0(r)-1}^{m_0(r)}  \frac{\ms d_l \big(\Phi^{-1}(t) \big)}{\ms d_\phi \big(\Phi^{-1}(t) \big)} \DD t \geq \frac 1r
\end{align*}
implying
\begin{align}
\label{K56}
m_0(r) \lesssim \Phi \Big(\Big[\frac{\ms d_\phi}{\ms d_l} \Big]^- \big(r \big)\Big).
\end{align}
Substituting back using \eqref{K19}, the estimate \eqref{K20} becomes
\begin{align}
\label{K57}
\begin{split}
 k (r) &\gtrsim r\int_{m_0(r)}^\infty \frac{\ms d_l \big(\Phi^{-1}(t) \big)}{\ms d_\phi \big(\Phi^{-1}(t) \big)} \DD t \gtrsim r \int_{\big[\frac{\ms d_\phi}{\ms d_l} \big]^- (r )}^\infty \ms d_l(x) \DD x.
\end{split}
\end{align}
Using \eqref{K68}, the theorem is proved.
\end{proof}

\subsection[{The case $\delta_l+\delta_{\phi}=2$}]{The case \boldmath{$\delta_l+\delta_{\phi}=2$}}

The method used in the proof of \Cref{K32} can also be applied in some exceptional cases. In the situation of \Cref{K29}, if \eqref{K30} is satisfied with $\delta_l+\delta_\phi=2$ and $\delta_l>1$, then \cite[Corollaries~4.7 and 5.2]{pruckner.reiffenstein.woracek:sinqB-arXiv} state that $\rho_H = \frac{1}{2}$. We show now that $\rho_H$ can in fact be larger than $\frac 12$ if $\delta_l=\delta_\phi=1$. A heuristic reason why this happens is that $H$ is close to limit point case when $\delta_l=1$.

\begin{proposition}
\label{K61}
Let $H=H_{l,\phi}$ be a Hamburger Hamiltonian. Suppose that $(\phi_j)_{j=1}^\infty$ is eventually monotone, and that
\begin{align*}
l_j \asymp j^{-1} (\log j)^{-\nu}, \qquad |\phi_{j+1}-\phi_j| \asymp j^{-1}
\end{align*}
for $j$ sufficiently large, where $\nu \in (1,2)$. Then
%\begin{alignat*}{3}
%	&\epsilon \leq 1 \quad && 
%	\Rightarrow
%		&&\quad  \max_{|z|=r} \log \|W_H(z)\| \asymp r^{\frac{1}{1+\epsilon}}; \\[1ex]		
%		&\epsilon > 1 \quad &&
%	\Rightarrow
%		&&\quad  \max_{|z|=r} \log \|W_H(z)\| \asymp r^{\frac 12} (\log r)^{\frac{1-\epsilon}{2}}
%\end{alignat*}
\begin{align*}
 \max_{|z|=r} \log \|W_H(z)\| \asymp r^{\frac{1}{\nu}}
\end{align*}
for $r$ sufficiently large. In particular, $\rho_H=\frac{1}{\nu}$.
\end{proposition}

\begin{proof}[Proof of \Cref{K61}]
Using the notation from the beginning of the section, we have
\begin{align*}
\ms d_l(t)=t^{-1} (\log t)^{-\nu}, \qquad \ms d_\phi (t)=t^{-1}.
\end{align*}
Further, let 
\begin{align*}
\ms c_l(t) \DE (\nu-1)\int_t^\infty \ms d_l(s) \DD s =(\log t)^{1-\nu}, \qquad \ms c_\phi \DE \ms c_l.
\end{align*}
%In the setting of subsection 5.2 of \cite{pruckner.reiffenstein.woracek:sinqB-arXiv}, taking $\ms c_\phi \DE \ms c_l$, we have
%\begin{align*}
%\delta_l &=\delta_\phi=1, & \delta&=2, \\
%\alpha_l &=\nu, \quad \alpha_\phi=0, & \alpha&=\nu, \\
%\gamma_l &=\gamma_\phi=0, & \gamma&=0, \\
%\beta_l &=\beta_\phi=\nu-1, & \beta&=\nu-1.
%\end{align*}
A calculation shows that we are in the situation of \cite[Example 5.9]{pruckner.reiffenstein.woracek:sinqB-arXiv}, and we obtain the upper bound
\[
 \max_{|z|=r} \log \|W_H(z)\| \lesssim r^{\frac 1\nu}.
\]
The proof of the lower bound is the same as that of \Cref{K32}, except for a few modifications. We have
\[
\Phi (t) = \int_1^{t} \ms d_\phi (s) \DD s = \log t
\]
and thus $\Phi^{-1}(x)=e^x$. Defining $I_j^+$ and $I_j^-$ as in the proof of \Cref{K32}, we have
\begin{align*}
&|I_j^+| \sim \Phi^{-1}(j+w)-\Phi^{-1}(j+u) =e^j \big( e^w-e^u \big), \\
&|I_j^-| \sim \Phi^{-1}(j-u')-\Phi^{-1}(j-w') =e^j \big( e^{-u'}-e^{-w'} \big)
\end{align*}
We will later choose $u,w,u',w'\in (0,1)$ such that $|I_j^+|<|I_j^-|$ for sufficiently large $j$. Proceeding as in Step 3, we arrive at \eqref{K13} and hence $w,w'$ need to be chosen such that $c' (w+w') \leq \frac{\pi}{2}$. One possibility is to use the ansatz $u=\frac w2$, $u'=\frac{w'}{2}$, $w' \DE \frac{\pi}{4} \min \{1, \frac{1}{c'} \}$, then choosing $w$ very small. Then, by our construction, \eqref{K14} holds. \\[1ex]
We proceed with Step 4, where no adjustments are necessary. Hence \eqref{K17} holds, leading to
\begin{align*}
&\sqrt{\det \Omega \big(x_{\lceil \Phi^{-1}( m)\rceil},x_{\lceil \Phi^{-1}( n) \rceil} \big)} \gtrsim  \sum_{j=m+1}^{n} \int_{\Phi^{-1}(j+u)}^{\Phi^{-1}(j+w)} \ms d_l(x) \DD x \\
&=\sum_{j=m+1}^{n} \int_{e^{j+u}}^{e^{j+w}} x^{-1} (\log x)^{-\nu} \DD x= \sum_{j=m+1}^{n} \int_{j+u}^{j+w} t^{-\nu} \DD t \\
&\asymp \sum_{j=m+1}^{n} \int_{j}^{j+1} t^{-\nu} \DD t \asymp \int_m^n t^{-\nu} \DD t.
\end{align*}
Steps 6 and 7 of the proof of \Cref{K32} can again be copied without changes, since
\begin{align*}
\frac{\ms d_l \big(\Phi^{-1}(t) \big)}{\ms d_\phi \big(\Phi^{-1}(t) \big)} = t^{-\nu}
\end{align*}
is regularly varying. By \eqref{K56}, we have $m_0(r) \lesssim r^{\frac{1}{\nu}}$, and thus the asserted lower bound is given by \eqref{K57}.
\end{proof}

\begin{remark}
\begin{itemize}
\item[]
\item[$\rhd$] If $\nu=2$ in \Cref{K61}, then $\rho_H=\frac 12$, and 
\begin{align*}
r^{\frac 12} \lesssim  \max_{|z|=r} \log \|W_H(z)\| \lesssim r^{\frac 12}\log r. 
\end{align*}
The proof of this estimate is the same as that of \Cref{K61}, except the upper bound taken from \cite[Example 5.9]{pruckner.reiffenstein.woracek:sinqB-arXiv} has a different form.

\item[$\rhd$] If $\nu>2$ , we also have $\rho_H=\frac 12$. All we can say about the precise growth is that
\[
r^{\frac 12} (\log r)^{-\frac{\nu}{2}} \lesssim \max_{|z|=r} \log  \|W_H(z)\| \lesssim r^{\frac 12}(\log r)^{1-\frac{\nu}{2}} . 
\]
In order to see this, use \cite[Corollary 5.2]{pruckner.reiffenstein.woracek:sinqB-arXiv} for the lower bound and the second row of \cite[Theorem 4.6]{pruckner.reiffenstein.woracek:sinqB-arXiv} for the upper bound.
\end{itemize}
\end{remark}

\section[{Growth of the Nevanlinna matrix of a Jacobi operator}]{Growth of the Nevanlinna matrix of a Jacobi operator}
\label{S4}

In this section we consider Jacobi parameters in a critical case, where $\frac{a_n}{b_n}$ converges to plus or minus two. We use two kinds of regularity conditions -- on the one hand, $\gamma$-temperedness (which is more general but complicated), and on the other hand power asymptotics (which are quite specific but more illustrative). Our results on Jacobi operators trace back to our results on canonical systems, but the change of parameters is a nontrivial task.

\subsection[{Rewriting a Jacobi matrix to a Hamburger Hamiltonian}]{Rewriting a Jacobi matrix to a Hamburger Hamiltonian}

Let $\ms J$ be a Jacobi matrix with parameters $a_n \in \bb R$, $b_n>0$.
 The minimal operator of $\ms J$ is defined as the closure in $\ell^2 (\bb N_0)$ of the operator $u \mapsto \ms Ju$ on the domain
\[
\mc D \DE \big\{u \in \ell^2 (\bb N_0) \DF \, u_k=0 \text{ for almost all } k \big\}.
\] 
Let us explain how to define a Hamburger Hamiltonian corresponding to $\ms J$ such that the minimal operator of $\ms J$ is unitarily equivalent to a certain symmetric operator associated with $H$. \newline

There are three degrees of freedom for choosing $H$: an initial condition $(\cos \psi,\sin \psi) \cdot y(0)=0$ (given by some $\psi \in [0,\pi)$), and the parameters $l_1,\phi_1$, where $\phi_1 \not\equiv \psi \mod \pi$. Our default choice is 
\[
\psi \DE 0, \qquad l_1 \DE 1, \qquad \phi_1 \DE \frac{\pi}{2}.
\]
With these parameters fixed, lengths and angles of a Hamburger Hamiltonian $H$ are uniquely determined by
\begin{align}
\label{K11}
a_0 &= \tan \phi_2, \\
\label{K44}
a_n &= - \frac{\sin (\phi_{n+2} - \phi_{n})}{l_{n+1} \sin (\phi_{n+2} -\phi_{n+1} ) \sin(\phi_{n+1} - \phi_{n})}, \qquad n \geq 1, \\
\label{K50}
b_n&=\frac{1}{\sqrt{l_{n+1}l_{n+2}}|\sin (\phi_{n+2}-\phi_{n+1})|}, \qquad n \geq 0.
\end{align}
These formulae can be found as (3.16) and (3.17) in \cite{kac:1999}, taking into account an index shift.
If $H$ is defined in this way, the minimal operator of $\ms J$ is unitarily equivalent to the closure of the operator whose graph is defined as the set of pairs $(f;g)$ satisfying\footnote{Some technical details are left out.} $f'=JHg$ with $(1,0) \cdot f(0)=0$ and $\supp f$ compact in $[0,L)$. 

The Nevanlinna matrix is denoted by
\[
W(z)=\begin{pmatrix}
A(z) & C(z) \\
B(z) & D(z)
\end{pmatrix},
\]
with $A,B,C,D$ defined in agreement with \cite{akhiezer:1961,berg.szwarc:2014}. The monodromy matrix of the Hamburger Hamiltonian defined above is related to these functions by
\begin{align*}
W_H(z)=\begin{pmatrix}
C(z) & A(z) \\
-D(z) & -B(z)
\end{pmatrix}
=
\begin{pmatrix}
A(z) & -C(z) \\
-B(z) & D(z)
\end{pmatrix}
J^T
\end{align*}
We will also need the fact that the lengths $l_k$ can be expressed in terms of the orthogonal polynomials $p_k,q_k$ of the first and second kind:
\begin{align}
\label{K93}
l_{k+1}=p_k(0)^2+q_k(0)^2, \qquad k \geq 0,
\end{align}
for our choice of $\psi,l_1,\phi_1$, cf. \cite[(3.22)]{kac:1999}.

\subsection[{\boldmath{$\gamma$}-tempered Jacobi parameters}]{\boldmath{$\gamma$}-tempered Jacobi parameters}

Our aim in this section is to prove \Cref{K10} as part of the more general \Cref{K78}, which describes the growth of the Nevanlinna matrix for Jacobi parameters $a_n,b_n$ that are perturbations of more regular parameters $\tilde{a}_n,\tilde{b}_n$. More precisely, $\tilde{a}_n,\tilde{b}_n$ should be $\gamma$-tempered in the sense of \cite{swiderski.trojan:2023}\footnote{We are not going to consider $N$-periodic modulations, i.e., we assume $N=1$.}. We state the definition below.

\begin{definition} 
\label{K77}
Consider Jacobi parameters $\tilde{a}_n,\tilde{b}_n$ with
\begin{align}
\label{K75}
&\lim_{n \to \infty} \tilde{b}_n=\infty, \\
\label{K82}
&\lim_{n \to \infty} \frac{\tilde{b}_{n-1}}{\tilde{b}_n}=1, \\
\label{K76}
&\lim_{n \to \infty} \frac{\tilde{a}_n}{\tilde{b}_n}=\mu
\end{align}
for some $\mu \in \{-2,2\}$. Further, let $(\gamma_n)_{n=0}^\infty$ be a sequence of positive numbers tending to infinity, satisfying
\begin{align}
\label{K83}
\lim_{n \to \infty} \sqrt{\gamma_{n+1}}-\sqrt{\gamma_n}=0.
\end{align}
Then $\tilde{a}_n,\tilde{b}_n$ are called \emph{$\gamma$-tempered} if, in addition to \eqref{K75}-\eqref{K76}, all of the following sequences are of bounded variation\footnote{$(x_n)_{n=0}^\infty$ is of bounded variation if $\sum_{n=0}^\infty |x_{n+1}-x_n|<\infty$.}:
\begin{align}
&\Big(\sqrt{\gamma_n}\Big(\frac{\tilde{b}_{n-1}}{\tilde{b}_n}-1 \Big) \Big)_{n=1}^\infty, \Big(\sqrt{\gamma_n}\Big(\frac{\tilde{a}_n}{\tilde{b}_n}-\mu \Big) \Big)_{n=1}^\infty, \Big(\frac{\gamma_n}{\tilde{b}_n} \Big)_{n=1}^\infty, \\
&\big(\sqrt{\gamma_n}-\sqrt{\gamma_{n-1}} \big)_{n=1}^\infty, \Big(\frac{1}{\sqrt{\gamma_n}} \Big)_{n=1}^\infty, \Big(\gamma_n \Big(1-\frac{\tilde{b}_{n-1}}{\tilde{b}_n}+\sgn (\mu) \big(\frac{\tilde{a}_n}{\tilde{b}_n}-\delta \big) \Big) \Big)_{n=1}^\infty
\end{align}

\end{definition}

Jacobi parameters $a_n,b_n$ with power asymptotics \eqref{K53} (assume that we are in the critical case) are not necessarily $\gamma$-tempered themselves, but they are perturbations of the parameters
\begin{align}
\label{K89}
\begin{split}
\tilde{b}_n &=n^{\beta} \Big(x_0+\frac{x_1}{n}\Big), \\
\tilde{a}_n &=n^{\beta} \Big(y_0+\frac{y_1}{n}\Big)=n^{\beta} \Big(\pm \frac{x_0}{2}+\frac{y_1}{n}\Big),
\end{split}
\end{align}
which, for $\beta>1$, are $\gamma$-tempered with $\gamma_n \DE n$. We then have
\begin{align}
\label{K84}
b_n=\tilde{b}_n \big(1+\tilde{\chi}_n \big), \qquad a_n=\tilde{a}_n \big(1+\tilde{\mu}_n \big)
\end{align}
for certain sequences $\tilde{\chi}_n,\tilde{\mu}_n$.

In the following theorem we determine the growth of the Nevanlinna matrix for a class of Jacobi parameters in the critical case. By the above reasoning, this class includes the case of power asymptotics \eqref{K53}.

\begin{theorem}
\label{K78}
Let $a_n,b_n$ be of the form \eqref{K84}, where $\tilde{a}_n,\tilde{b}_n$ are Jacobi parameters that are $\gamma$-tempered with $(\gamma_n)_{n=0}^\infty$ as in \Cref{K77}, and $\tilde{\chi}_n,\tilde{\mu}_n$ are such that
\begin{align}
\label{K85}
\sum_{n=0}^\infty \sqrt{\gamma_n}\big(|\tilde{\chi}_n|+|\tilde{\mu}_n| \big)<\infty.
\end{align}
 Assume that the limit
\begin{align*}
\tau_0 \DE \frac 14 \lim_{n \to \infty} \gamma_n \frac{\tilde{a}_n^2-4\tilde{b}_{n-1}\tilde{b}_n}{\tilde{b}_n^2},
\end{align*}
which exists by \cite[Theorem A]{swiderski.trojan:2023}, is negative.
Further assume
\begin{align*}
\sum_{n=0}^\infty \frac{\sqrt{\gamma_n}}{b_n}<\infty,
\end{align*}
and that $b_n$ and $\gamma_n$ are regularly varying with indices in $(0,2)$.
Then the Jacobi matrix $\ms J$ with parameters $a_n,b_n$ is in limit circle case and its Nevanlinna matrix $W$ satisfies
\[
	 \max_{|z|=r} \log\|W(z)\|  \asymp r \sum_{n=h(r)}^\infty \frac{\sqrt{\gamma_n}}{b_n},
	\]
	where $h(r) \DE \max \big\{n \in \bb N \cup \{0\} \DS \frac{b_n}{\gamma_n}<r\big\}$.
\end{theorem}

The following auxiliary lemma was communicated by G. Świderski.

\begin{lemma}
\label{K81}
Consider the situation of \Cref{K78}, except that we do not assume regular variation of $b_n$ and $\gamma_n$. Then there exists a solution $(f_n)_{n=0}^\infty$ of 
\begin{align}
\label{K74}
b_{n}u_{n+1}+a_nu_n+b_{n-1}u_{n-1}=0, \qquad n \geq 1,
\end{align}
such that 
\begin{align}
\label{K79}
|f_n|^2 \sim c \cdot \frac{\sqrt{\gamma_n}}{b_n}, \qquad n \to \infty,
\end{align} 
where $c>0$.
Furthermore, $(f_n)_{n=0}^\infty$ and $(\overline{f_n})_{n=0}^\infty$ are linearly independent and $\lim_{n \to \infty} (\arg f_n - \arg f_{n-1})$ exists modulo $\pi$ and is equal to $0$. 
\end{lemma}
\begin{proof}
Let us first prove the assertion in the unperturbed case where $\tilde{\chi}_n=\tilde{\mu}_n=0$ and thus $a_n,b_n$ are themselves $\gamma$-tempered. We explain the necessary adjustments for general $(\tilde{\chi}_n)_{n=0}^\infty$, $(\tilde{\mu}_n)_{n=0}^\infty$ satisfying \eqref{K85} in the end. 

We use ideas from \cite{swiderski.trojan:2023} and start similarly to Section 5 of that paper. We choose $K \DE \{0\}$ and identify all functions on $K$ appearing in \cite{swiderski.trojan:2023} with their value at $0$. Also note that $N=1$, $i=0$, and $\alpha_0=1$. Take $Y_j$ as in (5.1) of \cite{swiderski.trojan:2023}, which can be represented as
\begin{align}
\label{K86}
Y_j = \varepsilon \Big(I+\frac 1{\sqrt{\gamma_{j}}}R_j \Big)
\end{align}
where $\varepsilon=-\sgn \mu$ and $R_j$ converges to some $2 \times 2$-matrix $\mc R_0$ with negative discriminant. The matrix $Y_j$ has two eigenvalues $\lambda_j$ and $\overline{\lambda_j}$, where
\begin{align*}
\lambda_j=\varepsilon \Big(1+\frac{1}{\sqrt{\gamma_{j}}}\xi_j \Big)
\end{align*}
with $\xi_j,\overline{\xi_j}$ being the eigenvalues of $R_j$, chosen such that $\IM \xi_j>0$.  They converge to the eigenvalues $\xi, \overline{\xi}$ of $\mc R_0$, where $\IM \xi>0$.
 \\

Now we jump over to the proof of \cite[Theorem~4.1]{swiderski.trojan:2023}, whose reasoning we will follow except that we use $Y_j,R_j$ as defined above and\footnote{$\xi_j^-,\lambda_j^-$ are not needed here.} $\xi_j^+ \DE \xi_j,\, \lambda_j^+ \DE \lambda_j$. We obtain a sequence of vectors $(\phi_n)_{n=1}^\infty$ that satisfy 
\begin{align}
\label{K87}
\phi_{n+1} =\begin{pmatrix}
0 & 1 \\
-\frac{b_{n-1}}{b_n} & \frac{-a_n}{b_n}
\end{pmatrix} \phi_n, \qquad n>1,
\end{align}
and (4.8) of \cite{swiderski.trojan:2023}, i.e.,
\begin{align*}
\lim_{n \to \infty} \frac{\phi_n}{\prod_{k=n_0}^{n-1} \lambda_k} = p
\end{align*}
for some $p \neq 0$. Clearly $\phi_n$ is of the form $\binom{f_{n-1}}{f_{n}}$ where $(f_n)_{n=0}^\infty$ satisfies \eqref{K74}. Since $\lim_{n \to \infty} \lambda_n=\pm 1$ we see that $p$ is a nonzero multiple of $(1,1)^\top$ or of $(1,-1)^\top$, hence both entries of $p$ are nonzero. It follows from \cite[Claim~5.2]{swiderski.trojan:2023} that
\begin{align*}
|f_n|^2 \sim c \cdot \frac{\sqrt{\gamma_n}}{b_n}, \qquad n \to \infty,
\end{align*} 
with $c > 0$. Moreover, modulo $\pi$ the expression
\begin{align*}
d_n \DE \arg f_n-\sum_{k=n_0}^{n-1} \arg (\varepsilon\lambda_k )
\end{align*}
converges. Note that, since $\lim_{j \to \infty} \varepsilon\lambda_j=1$, for sufficiently large $j$ we get
\begin{align}
\label{K80}
\arg (\varepsilon \lambda_n)=\arctan \frac{\frac{\IM \xi_n}{\sqrt{\gamma_{n}}}}{1+\frac{\RE \xi_n}{\sqrt{\gamma_{n}}}} \sim \frac{\IM \xi}{\sqrt{\gamma_{n}}}.
\end{align}
This proves that 
\begin{align*}
&\arg f_n - \arg f_{n-1} \\
&= \Big(\arg f_n-\sum_{k=n_0}^{n-1} \arg (\varepsilon\lambda_k ) \Big)-\Big(\arg f_{n-1}-\sum_{k=n_0}^{n-2} \arg (\varepsilon\lambda_k ) \Big)+\arg(\varepsilon\lambda_{n-1})
\end{align*}
converges modulo $\pi$ to zero. To show that $(f_n)_{n=0}^\infty$ and $(\overline{f_n})_{n=0}^\infty$ are linearly independent we rewrite
\begin{align*}
\arg f_n - \arg \overline{f_n} \equiv 2\arg f_n \equiv 2\Big(d_n+\sum_{k=n_0}^{n-1} \arg (\varepsilon\lambda_k ) \Big) \mod \pi
\end{align*}
and argue that the right hand side is divergent. This holds because \eqref{K83} and the Stolz-Ces\`{a}ro theorem imply $\lim_{n \to \infty} \frac{\sqrt{\gamma_n}}{n}=0$, and by recalling \eqref{K80}. 

We still need to treat parameters $a_n,b_n$ of the form \eqref{K84}, where $\tilde{a}_n,\tilde{b}_n$ are $\gamma$-tempered with $(\gamma_n)_{n=0}^\infty$ and $\tilde{\chi}_n,\tilde{\mu}_n$ satisfy \eqref{K85}. Let $Y_j$ and $\tilde{Y}_j$ be defined by (5.1) of \cite{swiderski.trojan:2023} for $a_n,b_n$ and $\tilde{a}_n,\tilde{b}_n$, respectively.
Then, by \eqref{K86} and the formula above (10.3) in \cite{swiderski.trojan:2023} we can write
\begin{align*}
\tilde{Y}_j = \varepsilon \Big(I+\frac 1{\sqrt{\gamma_{j}}}\tilde{R}_j \Big), \qquad Y_j= \varepsilon \Big(I+\frac 1{\sqrt{\gamma_{j}}}\tilde{R}_j +V_j\Big),
\end{align*}
where $\tilde{R}_j$ converges to a matrix with negative discriminant and $\sum_{j=1}^\infty \|V_j\|<\infty$. For large enough $j$, let $\tilde{\xi}_j$ be the eigenvalue of $\tilde{R_j}$ whose imaginary part is positive and set
\[
\tilde{\lambda}_j \DE \varepsilon \Big(1+\frac{1}{\sqrt{\gamma_{j}}}\tilde{\xi}_j \Big).
\]
Since \cite[Theorem~4.4]{swiderski.trojan:2022} allows to perturb $\tilde{Y}_j$ with the summable sequence $\varepsilon V_j$, proceeding as in the proof of \cite[Theorem~4.1]{swiderski.trojan:2023}
we obtain a sequence $\phi_n$ satisfying \eqref{K87} for the perturbed parameters $a_n,b_n$ but with the same asymptotics as we would get for the solution of the unperturbed equation:
\begin{align}
\label{K88}
\lim_{n \to \infty} \frac{\phi_n}{\prod_{k=n_0}^{n-1} \tilde{\lambda}_k} = \tilde{p}.
\end{align}
Having \eqref{K88}, the rest of the proof is analogous to the unperturbed case.
\end{proof}

The solution $(f_n)_{n=0}^\infty$ of the Jacobi recurrence provided by \Cref{K81} gives us the information about the associated Hamburger Hamiltonian that is necessary to prove \Cref{K78}.

\begin{proof}[{Proof of \Cref{K78}}]

Our strategy is to apply \Cref{K29}, (iii), to the Hamburger Hamiltonian $H_{l,\phi}$ associated with $\ms J$. \vspace{8pt}

	{\it Step 1: Computing the decay of $l_n$ and $|\sin (\phi_{n+1}-\phi_n)|$.} Let $(f_n)_{n=0}^\infty$ be the solution of \eqref{K74} obtained from \Cref{K81}. Since $(\overline{f_n})_{n=0}^\infty$ is another solution of \eqref{K74} which is linearly independent of $(f_n)_{n=0}^\infty$, there is an invertible matrix $T$ such that $(p_n(0),q_n(0))=(f_n,\overline{f_n})\,T$ for all $n$, where $p_n,q_n$ are the orthogonal polynomials of the first and second kind associated with $\ms J$. By \eqref{K93} and \eqref{K79} we have 
\begin{align}
\label{K69}
l_{n+1}=\big\|(p_n(0),q_n(0)) \big\|_2^2 = \big\|(f_n,\overline{f_n})\, T \big\|_2^2 \asymp \big\| (f_n,\overline{f_n}) \big\|_2^2 \asymp \frac{\sqrt{\gamma_n}}{b_n}.
\end{align}
Relation \eqref{K50} thus implies
\begin{align}
\label{K70}
\sin (\phi_{n+1}-\phi_n) \asymp \frac{1}{\sqrt{\gamma_n}},
\end{align}
where the angles $\phi_n$ are normalised such that $\phi_{n+1}-\phi_n \in (0,\pi)$.
Having established \eqref{K69} and \eqref{K70}, it remains to prove that the angles $\phi_n$ have representatives $\hat{\phi}_n$ modulo $\pi$ that satisfy $|\hat{\phi}_{n+1}-\hat{\phi}_n| \asymp \frac{1}{\sqrt{\gamma_n}}$ and are again monotone.\vspace{8pt}
	
	{\it Step 2: Proving $\lim_{n \to \infty} \frac{l_{n+1}}{l_n}=1$.} By \Cref{K81} we also have
	\[
	\arg f_n \equiv \vartheta_n \mod \pi
	\]
	where $\vartheta_n$ is a sequence with $\lim_{n \to \infty} \vartheta_{n}-\vartheta_{n-1} =0$. From uniform continuity of the $\pi$-periodic function $\vartheta \mapsto \|(e^{i\vartheta},e^{-i\vartheta}) \, T \|_2$ we get
	\begin{align*}
	\lim_{n \to \infty} \bigg(\frac{l_{n+1}}{|f_n|^2}&-\frac{l_{n}}{|f_{n-1}|^2} \bigg)= \\
	&=\lim_{n \to \infty} \Big(\big\|(e^{i\vartheta_{n}},e^{-i\vartheta_{n}}) \, T \big\|_2 - \big\|(e^{i\vartheta_{n-1}},e^{-i\vartheta_{n-1}}) \, T \big\|_2 \Big)=0
	\end{align*}
	which, by multiplying with $\frac{|f_n|^2}{l_n} \asymp 1$, leads to 
	\begin{align}
	\label{K72}	
	\lim_{n \to \infty} \frac{l_{n+1}}{l_n}=\lim_{n \to \infty} \frac{|f_n|^2}{|f_{n-1}|^2}=1.
	\end{align}
	The last equality follows from \eqref{K79} since $\frac{b_{n-1}}{b_n} \sim 1$ by \eqref{K82}, $\frac{\gamma_{n}}{\gamma_{n-1}} \sim 1$ by \eqref{K83}.
	\vspace{2pt}
	
	{\it Step 3: Verifying monotonicity of angles.} By \eqref{K44},
	\[
	a_n=-b_n \sqrt{\frac{l_{n+2}}{l_{n+1}}} \frac{\sin (\phi_{n+2}-\phi_{n})}{\sin (\phi_{n+1}-\phi_{n})}=-b_{n-1} \sqrt{\frac{l_{n}}{l_{n+1}}} \frac{\sin (\phi_{n+2}-\phi_{n})}{\sin (\phi_{n+2}-\phi_{n+1})} .
	\]
	Using \eqref{K72} and the assumption $\lim_{n \to \infty} \frac{a_n}{b_n}=\lim_{n \to \infty} \frac{a_n}{b_{n-1}} =\mu$, we infer that
	\begin{align}
	\label{K73}
	\lim_{n \to \infty} \frac{\sin (\phi_{n+2}-\phi_{n})}{\sin (\phi_{n+1}-\phi_{n})}=\lim_{n \to \infty} \frac{\sin (\phi_{n+2}-\phi_{n})}{\sin (\phi_{n+2}-\phi_{n+1})} =-\mu \in \{-2,2\}.
	\end{align}
	Since $\phi_{n+1}-\phi_n \in (0,\pi)$ and
	\[
	\sin (\phi_{n+1}-\phi_{n-1})=\sin (\phi_{n+1}-\phi_{n})\cos  (\phi_{n}-\phi_{n-1})+\cos  (\phi_{n+1}-\phi_{n})\sin (\phi_{n}-\phi_{n-1}),
	\]
	relation \eqref{K73} implies that $\cos  (\phi_{n}-\phi_{n-1})$ and $\cos  (\phi_{n+1}-\phi_{n})$ have the same sign for all sufficiently large $n$. If this sign is positive, \eqref{K70} shows that 
	\[
	\phi_{n+1}-\phi_n \asymp \frac{1}{\sqrt{\gamma_n}}.
	\]
	In this case we set $\hat{\phi}_n \DE \phi_n$.
	
	If the sign of $\cos  (\phi_{n+1}-\phi_{n})$ is negative, then
	\[
	\pi-\big(\phi_{n+1}-\phi_n \big) \asymp \frac{1}{\sqrt{\gamma_n}}.
	\]
	Setting $\hat{\phi}_n \DE \phi_n-\pi n$, we have $H_{l,\phi}=H_{l,\hat{\phi}}$. Now $(\hat{\phi}_n)_{n=1}^\infty$ is decreasing, and we have
	\[
	- \big( \hat{\phi}_{n+1}-\hat{\phi}_n \big) \asymp \frac{1}{\sqrt{\gamma_n}}.
	\]\vspace{8pt}
	
	{\it Step 4: Applying \Cref{K29}.} In the previous steps we showed that item (iii) of \Cref{K29} is applicable with $\ms d_l(j)=\frac{\sqrt{\gamma_j}}{b_j}$ and $\ms d_\phi (j)=\frac{1}{\sqrt{\gamma_j}}$. A straightforward calculation shows that 
	\[
	 \max_{|z|=r} \log \|W(z)\| =  \max_{|z|=r} \log \|W_{H_{l,\phi}}(z)\| \asymp \ms m(r) \asymp r \sum_{n=h(r)}^\infty \frac{\sqrt{\gamma_n}}{b_n}
	\]
	which concludes the proof.
\end{proof}

\subsection[{Jacobi parameters with power asymptotics}]{Jacobi parameters with power asymptotics}

In this section we explain how to apply \Cref{K78} to Jacobi parameters $a_n,b_n$ having power asymptotics of the form \eqref{K53}. We assume that the simply critical case takes place, i.e., 
\[
\beta_2=\beta_1 \ED \beta \quad \wedge \quad |y_0|=2x_0 \quad  \wedge \quad \beta \neq \sigma \DE \frac{2x_1}{x_0}-\frac{2y_1}{y_0}.
\]
The aim of this section is to prove the following theorem.

\begin{theorem}
\label{K90}
Let $\ms J$ be the Jacobi matrix with parameters $a_n,b_n$ satisfying \eqref{K53}, assume that the simply critical case takes place and that $\sum_{n=1}^\infty \sqrt{n} \big(|\chi_n|+|\mu_n| \big)<\infty$. Then $\ms J$ is
\begin{itemize}
\item[-] in limit point case if $\beta>\sigma$, or $\beta<\sigma$ and $\beta \leq \frac 32$;
\item[-] in limit circle case if $\frac 32<\beta<\sigma$. In that case, for sufficiently large $r$ we have
\begin{align}
\label{K91}
 \max_{|z|=r} \log\|W(z)\| \asymp r^{\frac{1}{2 (\beta -1)}}
\end{align}
if $\frac 32<\beta<2$ and
\begin{align}
\label{K92}
 \max_{|z|=r} \log \|W(z)\| \asymp r^{\frac{1}{\beta}}
\end{align}
if $\beta > 2$. Furthermore, we have $\rho=\frac 12$ if $\beta=2$.
\end{itemize}
\end{theorem}

In order to determine when limit circle case takes place, we only need to apply results from \cite{swiderski.trojan:2023}, while for determining the growth of $W$ we need to show that \Cref{K78} is applicable.

Let the unperturbed parameters $\tilde{a}_n,\tilde{b}_n$ be given by \eqref{K89}, i.e.,
\begin{align*}
\begin{split}
\tilde{b}_n &=n^{\beta} \Big(x_0+\frac{x_1}{n}\Big), \\
\tilde{a}_n &=n^{\beta} \Big(y_0+\frac{y_1}{n}\Big)=n^{\beta} \Big(\pm \frac{x_0}{2}+\frac{y_1}{n}\Big).
\end{split}
\end{align*}
Then $a_n,b_n$ can be expressed by \eqref{K84} with $\tilde{\chi}_n=\chi_n \cdot \big(x_0+\frac{x_1}{n} \big)^{-1}$ and $\tilde{\mu}_n=\mu_n \cdot \big(y_0+\frac{y_1}{n} \big)^{-1}$, which satisfy
\[
\sum_{n=1}^\infty \sqrt{n} \big(|\tilde{\chi}_n|+|\tilde{\mu}_n| \big)<\infty.
\]

The following lemma is proved by somewhat lengthy but simple calculations, which we omit.

\begin{lemma}
For $\beta_1 > 1$, the parameters $\tilde{a}_n,\tilde{b}_n$ are $\gamma$-tempered with $\gamma_n \DE n$. With the notation from \Cref{K78} we then have $\tau_0 =\beta-\sigma$.
\end{lemma}

\begin{proof}[{Proof of \Cref{K90}}]
We first note that $\ms J$ is always in limit point case if $\beta \leq 1$, due to Carleman's condition. For $\beta>1$, since $a_n,b_n$ are perturbations of $\tilde{a}_n,\tilde{b}_n$ of the type considered in \cite[Section 10]{swiderski.trojan:2023}, it is enough to characterise limit circle/limit point case for the unperturbed parameters. It is easy to see that the function $\tau$ defined in \cite[(1.6)]{swiderski.trojan:2023} is constant equal to $\tau_0=\beta-\sigma$. If $\tau_0<0$, \cite[Theorem~9.1]{swiderski.trojan:2023} implies that $\ms J$ is in limit point case if and only if $\beta \leq \frac 32$, while for $\tau_0>0$ it is always in limit point case due to \cite[Theorem~9.2]{swiderski.trojan:2023} (to make the abstract statement more explicit, we can use (2.23) and (2.21) in the cited paper). Assume now $\frac 32<\beta<\min \{2,\sigma \}$, in which case limit circle case holds. To obtain \eqref{K91} we apply \Cref{K78}, noting that $h(r) \asymp r^{\frac{1}{\beta-1}}$ and hence
\begin{align*}
	 \max_{|z|=r} \log \|W(z)\|  \asymp r \sum_{n=h(r)}^\infty \frac{\sqrt{\gamma_n}}{b_n} \asymp r \int_{r^{\frac{1}{\beta-1}}}^\infty x^{\frac 12-\beta} \DD x \asymp r^{\frac{1}{2(\beta-1)}}.
\end{align*}
Finally, if $\beta \geq 2$ we note that \eqref{K69} and \eqref{K70} still hold and the associated Hamburger Hamiltonian satisfies \eqref{K30} with 
\[
\ms d_l(j)=\frac{\sqrt{\gamma_j}}{b_j} \sim j^{\frac 12-\beta}, \qquad \ms d_\phi(j)=\frac{1}{\sqrt{\gamma_j}} \sim j^{-\frac 12}.
\]
If $\beta=2$ then \cite[Corollaries~4.7 and 5.2]{pruckner.reiffenstein.woracek:sinqB-arXiv} imply $\rho=\frac 12$, while for $\beta>2$ we get the more precise formula \eqref{K92} from \cite[Theorem~5.3]{pruckner.reiffenstein.woracek:sinqB-arXiv}
\end{proof}

\subsection*{Acknowledgements.}
\noindent I am indebted to my supervisor Harald Woracek for his guidance and many fruitful conversations, to Grzegorz Świderski for patiently answering my questions, and to Per Alexandersson for helping me with some computations in the context of \Cref{K37}.

%---------
%   FINISH
%---------

\printbibliography

{\footnotesize

\begin{flushleft}
	J.~Reiffenstein \\
	Department of Mathematics\\
	Stockholms universitet \\
    106 91 Stockholm \\
	SWEDEN \\
	email: name dot surname at math.su.se \\[5mm]
\end{flushleft}

}

\end{document}